\documentclass[12pt]{article}
\usepackage[mathscr]{eucal}
\usepackage{amsbsy}
\usepackage{amssymb}
\usepackage{amsmath}
\usepackage{amsthm}
\usepackage{graphics}
\usepackage{epsfig}

\usepackage{latexsym}

\newcommand{\HH}{\mathrm{H}}

\newcommand{\PP}{\mathbf{P}}

\def\phi{\varphi}

\def\bbr{{\mathbb R}}
\def\bbb{{\mathbb B}}

\def\bz{{\mathbb Z}} 

\def\bll{{\bf l}}

\def\bk{{\bf k}}
\def\bt{{\bf t}}
\def\bs{{\bf s}}
\def\bll{{\bf l}}

\def\Bl1{{\bf 1}}
\def\B2{{\bf 2}}
\def\B0{{\bf 0}}

\def\x{\xi}
\def\e{\varepsilon}

\def\l{\lambda}

\def\s{\sigma}

\def\=A8{\"o}

\def\pas{\stackrel{a.s.}{\longrightarrow}}

\def\vr{\mathop{\hbox{\rm Var}}\nolimits}
\def\co{\mathop{\hbox{\rm conv}}\nolimits}
\def\card{\mathop{\hbox{\rm card}}\nolimits}

\newcommand{\beq}{\begin{equation}}
\newcommand{\eeq}{\end{equation}}
\newcommand\beqn{\begin{displaymath}}  
\newcommand\eeqn{\end{displaymath}}

\newcommand{\halmos}{\vspace{3mm} \hfill \mbox{$\Box$}\\[2mm]}
\theoremstyle{plain}
\newtheorem{teo}{Theorem}
\newtheorem{lem}[teo]{Lemma}

\newtheorem{cor}[teo]{Corollary}

\theoremstyle{definition}

\newtheorem{remark}[teo]{Remark}

\begin{document}

\title{On the asymptotic of convex hulls of Gaussian fields\footnote{The first named author was supported by grant VIZIT-3-TYR-013 of   Lithuanian
 Research Council.}
\footnotemark[0]\footnotetext[0]{ \textit{Short title:} asymptotic
of convex hulls of Gaussian fields }
\footnotemark[0]\footnotetext[0]{%
\textit{MSC 2000 subject classifications}. Primary 62G15, secondary
62G60, 60F15 .} \footnotemark[0] \footnotemark[0]\footnotetext[0]{ \textit{Corresponding
author:} Vygantas Paulauskas, Department of Mathematics and
 Informatics, Vilnius university, Naugarduko 24, Vilnius 03225, Lithuania,
 e-mail:vygantas.paulauskas@mif.vu.lt}
}


\author{ Youri Davydov$^{\text{\small 1}}$  and  Vygantas Paulauskas$^{\text{\small 2}}$ \\
{\small $^{\text{1}}$ Universit\'e Lille 1, Laboratoire Paul
Painlev\'e}\\
 {\small $^{\text{2}}$ Vilnius University, Department
of Mathematics and Informatics}\\
 {\small and Institute of Mathematics and Informatics,  }
 }



\maketitle

\begin{abstract}
We consider a Gaussian field $X = \{X_{\bt}, \;\bt \in T\}$  with
values in a Banach space  $\bbb$ defined on  a parametric
 set $T$ equal to $\bbr^m$ or $\bz^m.$
It is supposed that the  distribution $\cal P$ of $X_{\bt}$ is
independent of $t.$ We consider the asymptotic behavior of closed
convex hulls
$$
W_n = \co \{\,X_{\bt},\;\;\bt \in T_n\,\}
$$
where $(T_n)$ is an increasing sequence of subsets of $T$ and we
show that under some conditions of the weak dependence with
probability 1
$$
\lim_{n\rightarrow \infty} \frac{1}{b_n}\,W_n = {\cal E}
$$
(in the sense of Hausdorff distance), where the limit shape ${\cal E}$ is  the concentration ellipsoid of $\cal P.$

The asymptotic behavior of the mathematical expectations $Ef(W_n),$
where $f$ is an homogeneous function is also studied.

\end{abstract}
\vfill
\eject
\section{Introduction and formulation of results}


Let $\mathbb{B}$ be a separable Banach space and let $X =
\{X_{\bt},\;\; \bt \in T\,\}$ be a centered Gaussian process with
values in  $\mathbb{B}$ defined on some probability space
$\left\{\Omega, {\cal F}, \PP\right\}$. $T$ is some parametric
space, in our paper we shall consider two cases: $T = \bbr^m$ or $T
= \bz^m.$ For $\bt=(t_1, \dots, t_m), \ \bs=(s_1, \dots , s_m) \in
T$ (in both cases ) $|\bt-\bs|= \max_{1\le k\le m}|t_k-s_k|.$ In all
paper we shall assume that the marginal distributions of $X_\bt$ are
the same for all $\bt\in T$ and will be denoted by ${\cal P}$. The
measure ${\cal P}$ is Gaussian measure on a separable Banach space,
so by $\HH$ we denote the reproducing Hilbert space of this measure
and ${\cal E}$ will stand for the ellipsoid of concentration of the
measure ${\cal P}$ (i.e., the unit ball in $\HH$).

Let $(T_n)$ be an increasing sequence (this will be always
understood as $T_n\subset T_{n+1}$ ) of subsets of $T$ with $\nu_n
 \rightarrow \infty,$ where, in the case of $T =
\bz^m$, \ $\nu_n$ is defined as  $ \card \{T_n \}$, while in the
case of $T = \bbr^m$, \ $\nu_n = \l^m(T_n),$  where $\l^m$ denotes
the Lebesgue measure in $\bbr^m$. For a set $A\subset \mathbb{B}$
let us denote by $\co \{A\}$ the closed convex hull of the set $A$.
We are interested in the limit  behavior of the sequence of sets
$$
W_n = \co \{\,X_{\bt},\;\;\bt \in T_n\}.
$$
This problem is interesting and important since it can be
considered as the multivariate generalization of classical and
deeply investigated problem on the behavior of extreme values of
Gaussian processes, see, for example, \cite{Leadbetter},
\cite{Berman}, \cite{Mittal} and references there. The limit
behavior of sets $W_n$ is closely related with the limit behavior
of Gaussian samples, see \cite{Goodman}, and has various
interesting applications, see \cite{Randon}.

In
\cite{Davydov} the case $T=\bz$ and $X=\{X_1, X_2, \dots\}$ being
independent and identically distributed (i.i.d.) random variables
with values in $\bbb$ was studied, while in \cite{DavDom} the case
of stationary sequences with
$T = \bbr$ and $\bbb=\bbr^d$ was considered. It was stated in particular that under mild conditions
with probability 1
$$
\lim_{n\rightarrow \infty} \frac{1}{\sqrt{2\ln n}}\,W_n = {\cal E}
$$
(in the sense of Hausdorff distance), where $T_n = T\cap [0, n]$  and the limit shape ${\cal E}$ is
the concentration ellipsoid
 defined by the covariance structure of $X$.
We generalize and complement the
statements of \cite{Davydov} and \cite{DavDom}.

 In order to formulate our results we need some more notation.
 $\bbb^\ast$ will stand for the conjugate space of $\bbb$ and $\langle \cdot, \cdot\rangle$
 denote the bilinear form
defining the duality between $\bbb$ and its conjugate space.
$B_r(x)$ and $S_r(x)$ denote the  closed ball and the  sphere,
respectively, with radius $r$ and center $x\in \bbb$, while
$B_r^*(x)$ and $S_r^*(x)$ stand for corresponding sets in $\bbb^*.$

Since in our setting the sets $W_n$ will be compact, we introduce
the separable complete metric space $\mathcal{K}_{\mathbb B}$ of all
nonempty compact subsets of a Banach space $\mathbb B$ equipped with
the Hausdorff distance $\rho_{\mathbb B}$ :
$$
\rho_{\mathbb B}(A,B) = \max\{\inf\{\,\epsilon \;|\;A \subset
B^\epsilon\},\;\; \inf\{\,\epsilon \;|\;B \;\subset
A^\epsilon\}\},
$$
$A^\epsilon$ is the open $\epsilon$-neighbourhood of $A$.
Convergence of compact sets in $\bbb$ always will be in this
metric.

Also in all paper we use the notation $b(t) = \sqrt{2\ln
(t\vee2)}$. Our first result is in the case $T = \bz^m.$

\begin{teo}\label{thm1}
Suppose that a Gaussian process $X$ with  the same marginal
distributions ${\cal P}$ for all $t$ satisfies the following
condition of the weak dependence
\begin{equation}\label{condVD}
 \;\;\;\;\; \forall x^\ast \in \mathbb{B}^\ast\;\;\;
E\left\langle X_\bt,\,x^\ast\right\rangle\left\langle
X_\bs,\,x^\ast\right\rangle \rightarrow 0, \;\;\; as \ \ |\bt-\bs|
\rightarrow \infty.
\end{equation}
Then
\begin{equation}\label{convW}
\frac{1}{b(\nu_n)}W_n \stackrel{a.s.}{\rightarrow} {\cal E}, \quad
as \ n\to \infty,
\end{equation}
where $\pas $ denotes the convergence a.s. (and, as it was
mentioned, in the metric $\rho_{\bbb}$).
\end{teo}

In the continuous case ($T = \bbr^m$) we need two additional
hypothesis. Now it is not sufficient to require that marginal
distributions of the process $X$ are the same, and we  suppose that
our process is stationary. For the subsets
  $T_n$ (in discrete case it was finite sets) we assume that they are compact sets satisfying the condition
\begin{equation}\label{condVH}
 \;\;\;\;\; \forall \e >0 \;\;\;\lim_n\frac{\l^m((\partial
T_n)^\e)}{\l^m(T_n)} \rightarrow 0,
\end{equation}
where $\partial T_n$ stands for the boundary of $T_n.$

\begin{teo}\label{thm2}
Suppose that the process $\{X_\bt,\; \bt\in \bbr^m\}$ is stationary and
 the conditions (\ref{condVD}) and (\ref{condVH}) are
fulfilled. Then  $W_n\in \mathcal{K}_{\mathbb B}$ a.s. and the
relation (\ref{convW}) takes place.
\end{teo}

Having the results on the convergence a.s. (and, therefore, in
distribution),  we can easily obtain, as in \cite{Davydov},  the
convergence of mean values for various functionals of these distributions.

 \vspace{7pt}


Let $f: \mathcal{K}_\bbb \rightarrow \mathbb{R}$ be a continuous
non-negative increasing homogeneous function of degree $p$ , that
is \vspace{5pt}




$f(A) \geq 0 \;\;\; \forall A \in \mathcal{K}_\bbb;$
\vspace{5pt}

$f(A_1) \leq f(A_2) \;\;\;\forall A_1\subset A_2,\;\, A_1,A_2 \in \mathcal{K}_\bbb;$
\vspace{5pt}

$f(cA) = c^pf(A),\;\; \forall \;c\geq 0,\; \forall A \in \mathcal{K}_\bbb.$ \vspace{5pt}


\begin{teo} \label{esper}
 Let $f$ be a homogeneous function of degree $p$ with the properties described above. Suppose additionally that there
 exists a constant $C$ such that for all $A \in {\cal K}_\bbb$
 \beq\label{moment}
 f(A) \leq C[d(A)]^p,
 \eeq
 where $d(A) = \sup_{x,y \in A} \left\|x-y\right\|$ is the diameter of $A.$
 Then, under hypothesis of Theorem 1 or 2  for all $a >0$
$$
E\exp\left\{af^{\frac{2}{p}}\left(\frac{1}{b(\nu_n)}
\,W_n\right)\right\}\; \longrightarrow \;\;\exp\left\{af^{\frac{2}{p}}\left(W\right)\right\}.
$$
\end{teo}
\begin{cor} \label{cor1}
Let $f$ be a function with the pro\-per\-ties described above. Then, under hypothesis of Theorem 1 or 2 for all $m>0$
$$
Ef^m\left(\frac{1}{b(\nu_n)}\,W_n\right) \rightarrow f^m\left(W\right).
$$
\end{cor}


 This theorem and corollary give in particular the asymptotic behavior for mean values of all reasonable geometrical
 characteristics of $W_n$ (such as diameter or  volume and  surface measure in the case of finite-dimensional $\bbb$).


\section{Auxiliary lemmas }

The first lemmas are about compact sets in $\bbb$.

\vspace{5pt}

\begin{lem}\label{lem1} If $A \in {\cal K}_\bbb,$ then $\co (A) \in {\cal K}_\bbb$ and
the mapping\\ $\co : {\cal K}_\bbb \longrightarrow {\cal
K}_\bbb,\;\;\; A \rightarrow \co (A),\;$ is 1-Lipshitz : $$
\rho_{\bbb} (\co (A), \co (B)) \leq \rho_{\bbb} (A, B). $$
\end{lem}

The proof of the lemma is elementary and is left for a reader.
\vspace{5pt}

\begin{lem}\label{lem2} Suppose that $A_n, A \in {\cal K}_\bbb$ are such that for some
sequence $(\e_n)\downarrow 0$
$$
A_n \subset A^{\e_n},\;\;\; \forall n.
$$
Then $(A_n)$ is relatively compact in ${\cal K}_\bbb.$
\end{lem}

\vspace{5pt}

{\bf Proof.} Let $S_n=\{x_{n,j}, \ j=1, \dots, k_n \}$ be a
$\e_n$-net for the set $A$, then it will be $2\e_n$-net for the set
$A^{\e_n}$.  Let us consider the set
$$
{\tilde A}_n=(A\cap A_n)\cup \left \{x_{n,j}: B_{2\e_n}(x_{n,j})\cap
A\ne \emptyset \right \}.
$$
 From the construction of the sets ${\tilde A}_n$ we have that ${\tilde
 A}_n$ are compact sets and ${\tilde A}_n\subset A$ for all $n$. It
 is known that if a sequence of compact sets is inside of one fixed
 compact sets, then this sequence is relatively compact
 (see \cite{Schneider}, Th.1.8.4., for finite dimensional case; for Banach spaces
 we have no relevant reference, but the proof is analogous). Again, from the construction of the sets ${\tilde A}_n$ we
 have the following relations
 $$
A_n \subset {\tilde A}_n^{2\e_n}, \quad {\tilde A}_n \subset
A^{2\e_n},
 $$
whence it follows that $\rho_{\mathbb B}(A_n, {\tilde A}_n)\le
2\e_n$. Since the sequence ${\tilde A}_n$ is relatively  compact,
the same property has the  sequence $A_n,$ too. The lemma is proved. \halmos

\begin{lem}\label{lem3} If $(A_n),\;\; A_n \in {\cal K}_\bbb,$ is relatively
compact, then  the sequence $\{\co (A_n)\}$ is relatively compact,
too.
\end{lem}
This fact follows directly from Lemmas \ref{lem1} and \ref{lem2}.

\vspace{5pt}

\begin{lem}\label{lem4}  If $A = \cap_n A_n,$ where $A_n \in {\cal K}_\bbb$, and
$(A_n)\downarrow,$ then $A_n \rightarrow A.$
\end{lem}
\vspace{5pt}

{\bf Proof.} Let's assume the opposite. Then for some $\delta>0 $
there exists a subsequence $(n')$, for which $\rho_{\bbb}(A,
A_{n'})> \delta \quad \forall n'.$  Without restriction of a
generality we can suppose that $(n') = (n).$ Then, as $A_n$ can not
be a subset of $ A^\delta, $ we can find $x_n \in A_n $ such that
$d(x_n, A)> \delta$ (here $d(x, A)$ stands for a distance from a
point $x$ and a set $A$ in a Banach space $\bbb$). As $A_n \subset
A_1, $ and $A_1 $ is compact, it is possible to choose a subsequence
$(n_k) $, for which $x_{n_k} \rightarrow x_0. $ It is clear that for
the limit point $x_0$ we will have $d(x_0,A)\geq \delta.$ On the
other hand, for each $m$ and for all sufficiently large $k,$ \
$x_{n_k} \in A_m$, which means that $x_0 \in A_m, \; \forall m.$
Hence, $x_0$ must belong to $A$ in contradiction to the previous
conclusion.\halmos



\begin{lem}\label{L1}
Under conditions of Theorem \ref{thm1} the sequence
$\left\{\frac{1}{b(\nu_n)}W_n\right\}$ is relatively compact a.s.
\end{lem}

{\bf Proof.} Let's show that with probability 1 compact sets
$$
K_n = \left\{\frac{1}{b(\nu_n)}X_\bk, \bk \in T_n\right\}
$$
 form a relatively compact sequence in ${\cal
K}_\bbb $. Then, due to Lemma \ref{lem4}, we get the result.

Let us renumber r.v. $X_\bk $ with the  indices $\bk $ from $\cup_nT_n $
as follows: at first somehow (but in a row) let's enumerate the
random variables with indices lying in $T_1 $ (there will be $\nu_1
$ of them), then will add the indices corresponding to random
variables from $T_2\setminus T_1 $, and so on. The sequence obtained
in this way we will denote by $\{Z_n\}$.

As r.v. $Z_k $ have the same distribution, it is possible to use the
first part of Theorem 1 from \cite{Goodman} (in its proof the
assumption of independence isn't used),  which gives a.s.
convergence
$$
\max_{k\leq n}d(Z_k, b(n) {\cal E}) \rightarrow 0,\quad n
\rightarrow \infty.
$$
It means that
$$
K_n \subset {\cal E}^{2\varepsilon_n},
$$
where a.s.
$$
\varepsilon_n = \max_{k\leq b(\nu_n)}\left\{ d\left(
\frac{Z_k}{b(\nu_n)}, \; {\cal E}\right) \right\}  \rightarrow 0.
$$
As ${\cal E}$ is compact, we conclude the proof applying Lemma
\ref{lem2}.\halmos

\begin{lem}\label{L2}
Let $(\x_n)$ be a real-valued  Gaussian centered sequence with
$\vr(\x_n) = \sigma^2 \;\;\;\forall n.$ 
Let
$$
c= \liminf_n \left\{\frac{1}{b(n)}\max_{k\leq
n}\left\{\x_k\right\}\right\}.
$$
Suppose that  
$$
r = \sup_{n\neq
l}\frac{|E\x_n\x_l|}{\sigma^2} < 1/2.
$$
Then
$$
\sigma(\sqrt{1-r}-\sqrt{r}) \leq c \leq \sigma.
$$
\end{lem}



 {\bf Proof.} The upper bound $c
\leq \sigma$ is the well-known fact (see i.e. Lemma \ref{propD} below), and for the proof of the lower
bound we introduce independent standard Gaussian random variables
$\eta$ and $\zeta_k, k\ge 1$ and define
$$
{\tilde \x}_n=\s {\sqrt {1-r}}\zeta_n +\s {\sqrt {r}}\eta.
$$
Then
$$
\vr({\tilde \x}_n) = \sigma^2 \quad \mbox {and}\quad E{\tilde
\x}_n{\tilde \x}_m=\s^2 r\ge E\x_n \x_m, \ \forall n, m.
$$
Therefore, from Slepian lemma (see Corollary 3.12 in \cite{Ledoux})
it follows that
$$
\PP \left \{\max_{k\le n}\x_k \le \l\right \}\le \PP \left
\{\max_{k\le n}{\tilde \x}_k \le \l\right \}, \quad \forall \l.
$$
Denoting
$$
Z_n=\frac{1}{b_n}\max_{k\le n}\x_k, \quad {\tilde
Z}_n=\frac{1}{b_n}\max_{k\le n}{\tilde \x}_k,
$$
and taking $\l=\s s b_n$ with $s<{\sqrt {1-r}}-\sqrt{r}$  and
$b_n=b(n)$, we have \beq \label{Slep} \PP \left \{Z_n \le \s s
\right \}\le \PP \left \{{\tilde Z}_n \le \s s \right \}. \eeq It
remains to prove \beq \label{Slep1} \sum_n \PP\left \{ Z_n \le \s s
\right \}<\infty, \eeq since then by Borel-Cantelli lemma it will
follow that
$$
{\liminf}_n Z_n \ge \s s \quad \mbox {a.s.}
$$
Taking into account  the relation (\ref{Slep}) it is sufficient to
prove that \beq \label{Slep1a} \sum_n \PP\left \{{\tilde Z}_n \le \s
s \right \}<\infty, \ \mbox {for all} \ s<{\sqrt {1-r}}-\sqrt{r}.
\eeq For this aim we must to show that \beq \label{Slep2} \sum_n
J_n(s)<\infty, \ \mbox {for all} \ s<{\sqrt {1-r}}-\sqrt{r}, \eeq
where
$$
J_n\equiv J_n(s)=\int_R \Phi^n \left (\frac{sb_n-{\sqrt r}t}{\sqrt
{1-r}} \right )\phi(t)dt.
$$
Here $\phi$ and $\Phi$ are the density function  and distribution
function, respectively,  of a standard normal random variable. To
simplify the notation, we denote
$$
d=\frac{s}{\sqrt {1-r}}<1, \quad a=\sqrt {\frac{r}{1-r}}, \quad
$$
Then we can write
$$
J_n=\int_R \Phi^n \left (db_n-at\right )\phi(t)dt= a^{-1}\int_R
\Phi^n \left (db_n-y\right )\phi(y/a)dy.
$$

 Let us take a positive real number
$\e,$ which will be chosen later and write
$$
J_n=I_{1,n}+I_{2,n},
$$
where $I_{1,n}$ and $I_{2,n}$  are corresponding integrals over the
intervals $(-\infty, -\e b_n)$ and $(-\e b_n, \infty).$ In the first
interval we simply estimate $\Phi^n \left (db_n-y\right )\le 1$ and
we get
$$
I_{1,n}\le \int^{\infty}_{\e b_na^{-1}}\phi (t)dt\le \frac{Ca}{\e
b_n}\exp (-\e^2 b_n^2a^{-2}/2)=\frac{C}{\e_2 b_n} n^{-\e^2a^{-2}}.
$$
Here and in what follows $C$ stands for an absolute constant, not
necessary the same in different places. If we chose $\e$ satisfying
condition
\begin{equation}\label{condep1}
\e>a=\sqrt {\frac{r}{1-r}}
\end{equation}
then we get \beq \label{Slep3} \sum_n I_{1,n}<\infty. \eeq Let us
note that $\Phi^n \left (db_n-y\right )$ is the decreasing function
of $y$ , therefore \beq \label{Slep4} I_{2,n}\le \Phi^n \left
((d+\e)b_n\right ).\eeq  We have
$$
\Phi^n((d+\e)b_n)=\left (1-\int_{(d+\e)b_n}^\infty \phi (t)dt \right
)^n.
$$
Since $1-\Phi (z)\sim z^{-1}\phi(z)$ for $z\to \infty,$ there exists
a constant $c_1>0$ such that for sufficiently large $z,$\  $1-\Phi
(z)> c_1z^{-1}\phi(z)$, therefore, for sufficiently large $n$
\begin{eqnarray} \label{Slep5}
\Phi^n \left ((d+\e)b_n\right ) &\le & \left (1-\left (1-\Phi \left
((d+\e)b_n\right )\right )
\right )^n  \nonumber \\
 & \le & \exp \left \{-n \left (1-\left (1-\Phi
\left ((d+\e)b_n\right )\right )\right )\right \}\nonumber \\
& \le & \exp \left (-\frac{c_1n^{1-(d+\e)^2}}{(d+\e)b_n} \right ).
\end{eqnarray}

From (\ref{Slep4}) and (\ref{Slep5}) we obtain
$$
I_{2,n}\le \Phi^n \left ((d+\e)b_n\right )\le \exp \left
(-\frac{c_1n^{1-(d+\e)^2}}{(d+\e)b_n} \right ).
$$
Now, if we chose $\e$ satisfying condition
\begin{equation}\label{condep2}
\e<1-d=1- \frac{s}{\sqrt {1-r}},
\end{equation}
 then \beq \label{Slep6} \sum_nI_{2,n}\le \sum_n \exp \left
(-\frac{c_1n^{1-(d+\e)^2}}{(d+\e)b_n} \right )<\infty. \eeq
 It remains to note that due to the condition $s<{\sqrt
 {1-r}}-\sqrt{r}$  it is possible to choose $\e$,
 satisfying both conditions (\ref{condep1}) and (\ref{condep2}),
 since
 $$
\sqrt {\frac{r}{1-r}}< 1- \frac{s}{\sqrt {1-r}}.
 $$

Estimates (\ref{Slep3}) and (\ref{Slep6}) prove (\ref{Slep2}). The
lemma is proved. \halmos

\begin{remark}\label{remark}
At first we were sure that only simple estimates which we had used
 do not allow to prove stronger statement, namely,
under condition that $r<1$
\begin{equation}\label{strongstat}
\sigma\sqrt{1-r} \leq c \leq \sigma.
\end{equation}
It turned out that even exact investigation of the integrand
function
$$
\Phi^n \left (\frac{sb_n-{\sqrt r}}{\sqrt {1-r}} \right )\phi(t)
$$
does not allow to achieve this goal. Contrary, it is possible to
show (we do not provide these calculations since they are rather
lengthy) that for  $\sqrt {1-r}-\sqrt{r}\le s < {\sqrt {1-r}}$ the
series in (\ref{Slep2}) diverges. But since the divergence of this
series does not imply the divergence of series in (\ref{Slep1}), the
question if the above stated strengthening (\ref{strongstat}) of the
lemma is possible remains open.
\end{remark}

\begin{lem}\label{L4}
Let $(\x_\bk),\, \bk \in \bz^m,$ be a real-valued Gaussian centered
field with $\vr(\x_\bk) = \sigma^2 \;\;\;\forall\; \bk$ and
$r_{\bk,\bll} = E\x_\bk\x_\bll \rightarrow 0 \;$ as $|\bk-\bll|
\rightarrow \infty.$ Let $(T_n)$ be an increasing sequence of
subsets of $\bz^m$ with $\nu_n = \card \{T_n \} \rightarrow
\infty.$

Then
\begin{equation}\label{el2}
Z_n = \frac{1}{b(\nu_n)}\max_{\bk \in T_n}\left\{\x_\bk\right\}
\pas \sigma.
\end{equation}
\end{lem}

{\bf Proof.} Fix $\varepsilon \in (0,1/2).$ By condition there
exists $a>0$ such that $|r_{\bk,\bll}| < \varepsilon\sigma^2$ if
$|\bk-\bll| \geq a.$ We will show that it is possible to find an
increasing sequence $(\tilde{T}_n)$ of subsets of $\bz^m$ with the
following properties:

\begin{enumerate}

\item $\tilde{T_n} \subset T_n;$

\item $T_n \subset (\tilde{T}_n)^a;$

\item $\forall \; \bk,\bll \in \cup_n\tilde{T}_n, \; \bk \neq\bll,$ we have
$|\bk-\bll| \geq a.$
\end{enumerate}
From 1.--3. it follows that
$$
\tilde{\nu}_n:= \card\{\tilde{T_n} \} \leq \nu_n \leq
(2a)^m\tilde{\nu_n}.
$$
Therefore, $b(\tilde{\nu_n}) \sim b(\nu_n)$, and we have a.s.
$$
\liminf_n Z_n \geq \liminf_n \frac{1}{b(\nu_n)}\max_{\bk \in
\tilde{T_n}}\left\{\x_\bk\right\} \geq \sigma \varphi(\varepsilon),
$$
where $\varphi(r)= \sqrt{1-r}-\sqrt{r},$ by Lemma \ref{L2}.

As $\varphi(r) \rightarrow 1$ when $r \rightarrow 0,$ we deduce that
a.s. $\liminf_n Z_n \geq \sigma.$ The opposite inequality $\limsup_n
Z_n \leq \sigma$ being well known, we arrive to (\ref{el2}).

Now we provide the construction of the sequence $(\tilde{T}_n)$. For
a finite subset $B$ of $\bz^m$ denote by ${\cal L}_a(B)$ the family
$\{E,\;E\subset B\}$ of all subsets of $B$ such that $\forall \;
\bk,\bll \in E, \bk \neq\bll,$ we have $|\bk-\bll| \geq a.$ If
${\cal L}_a(B)$ is not empty, let $B_a$ be one of  its elements of
maximal cardinality. If ${\cal L}_a(B)$ is empty, we use the
notation $B_a$ for arbitrary chosen singleton $\{\bk\}\subset B.$ In any
case it is clear that
$$
B_a \subset B \subset (B_a)^a,
$$
which gives the inequalities
$$
\card(B_a) \leq \card(B) \leq (2a)^m\card(B_a).
$$
    We define our sequence by induction. We set $\tilde{T_1}= (T_1)_a.$   When $\tilde{T_n}$ is defined, then $\tilde{T}_{n+1}$
 is equal to $\tilde{T_n}$, if $T_{n+1} \subset (T_n)^a,$ and $\tilde{T}_{n+1}$ is equal to\\ $\tilde{T_n} \cup (T_{n+1}\setminus (T_n)^a)_a$
 in the case  when
$T_{n+1}\setminus (T_n)^a \neq \emptyset.$

It is easy to see that the properties 1.--3. are fulfilled.
Therefore the lemma is proved. \halmos

In the sequel we shall need the notion of a support function. The function
$\mathcal{M}_A(\theta),\;\;\theta \in S^*_1(0),$   defined by the
relation
$$
\mathcal{M}_A(\theta)= \sup_{x\in A}\langle x, \theta\rangle,
\;\;\;\theta \in S^*_1(0),
$$
is called a {\it support
function} of a set $A\in \mathcal{K}^d.$

 A compact convex set $A$ is characterized by its support function
since
\[
A=\bigcap_{\theta\in S^{\ast}(0,1)} \{u\in\bbb; \langle
u,\theta\rangle\leq \mathcal{M}_A(\theta)\}.
\]

It follows easily from definition that $\mathcal{M}_A$ is 1-Lipshitz
and that
$$
\rho_{\bbb}(A,B) = \sup_{\left\|
\theta\right\|=1}|\mathcal{M}_A(\theta) - \mathcal{M}_B(\theta)|.
$$

\begin{lem}\label{L5}
Let $(B_n)_{n\geq 0}$ be a sequence of random convex elements in ${\cal
K}_\bbb$. Assume that $(B_n)$ is a.s. relatively compact. Assume
also that there exists a (deterministic) function
$\phi:S_1^*(0)\to\bbr$ such that, for all $\theta \in
S_1^{\ast}(0)$,
\[
\mathcal{M}_{B_n}(\theta)\to \varphi(\theta)\quad  \mbox{a.s., as }
n\rightarrow +\infty.
\]
Then $\varphi$ is the support function of a
 set $A\in{\cal K}_\bbb$  and
\[
B_n\rightarrow A  \quad  \mbox{a.s., as } n\rightarrow +\infty.
\]
\end{lem}

{\bf Proof.} Let $\Omega_1, {\mathbb P}(\Omega_1)=1,$ be a subset of
$\omega$ for which the sequence  $(B_n)$ is relatively compact. Let
$D$ be a countable dense subset of $S^{\ast}_1(0)$ and $\Omega_2,
{\mathbb P}(\Omega_2)=1,$ be a subset of $\omega$ for which
$\mathcal{M}_{B_n}(\theta)\to \phi(\theta).$ Fix $\omega$ from
$\Omega_1\cap \Omega_2.$ Let $A\in{\cal K}_\bbb$ be a limit point of
the sequence $(B_n)_{n\geq 1}$.  We denote by $(m_n)_{n\geq 1}$ an
increasing sequence such that $B_{m_n}\rightarrow A.$ Then, for all
$\theta\in D$, $\mathcal{M}_{B_{m_n}}(\theta)\rightarrow
\mathcal{M}_A(\theta)$, as $n\rightarrow+\infty$. At the same time
$\mathcal{M}_{B_{m_n}}(\theta)\rightarrow \varphi(\theta).$
 Using uniqueness of the limit, we obtain the  equality $\mathcal{M}_A=\varphi$ a.s. on $D$.
 If $A'$ is another limit point of the sequence $(B_n)_{n\geq 1}$, we  have also $\mathcal{M}_{A'}=\varphi$  on $D$.
  Consequently $\mathcal{M}_{A}= \mathcal{M}_{A'}$ on $D$ and
 by continuity, the equality holds on  $S^{\ast}_1(0)$.
Finally, $(B_n)$ has a unique limit point and  $B_n\rightarrow A$ in
${\cal K}_\bbb$ almost surely as $n\rightarrow+\infty$. Since for
each $\omega \in \Omega_1\cap \Omega_2$ we get the same
deterministic support  function of the set $A$ we get
 that this limit point $A$ is  deterministic and its
support function $\mathcal{M}_A=\phi$. \halmos \vspace{10pt}

We will  need also the following
general result,
dealing with the maximum of sub-Gaussian random variables (see i.e. \cite{Davydov},  Lemmas 1 and 3 therein).
\begin{lem}\label{propD}
Let $(Y_n)_{n\geq 0}$ be a sequence of identically distributed
random variables such that for some $\zeta>0$
\[
 E[e^{\gamma Y_0^2}]<\infty, \quad \mbox{for all}\ \gamma< \frac{1}{2\zeta^2}.
\]
Let
$$
Z_n = \frac{1}{\sqrt{2\ln n}}\max\{\,Y_1,\ldots,Y_n\}.
$$
Then,

\[
 \limsup_{n\to+\infty} Z_n \leq \zeta \quad
 \mbox{a.s.},
\]

and  for any $a>0$,
$$
\limsup_n E\exp\left\{a Z_n^2\right\} < \infty.
$$
\end{lem}
Note that Lemmas 1 and 3 in \cite{Davydov} are stated for
independent random variables, but it is clear from the proof that
the assumption of independence is unnecessary.

\section{Proofs of Theorems}

{\bf Proof of Theorem \ref{thm1}.} It is easy to see that
$$
\mathcal{M}_{\cal E}(\theta) = \sqrt{E\left\langle X_\bk,
\theta\right\rangle^2},\;\; \theta \in S_1^{\ast}(0).
$$
Due to Lemmas \ref{L1} and \ref{L5} it is sufficient to show that
$\forall \theta \in S^*_1(0)$  \beq\label{proofTh1}
\mathcal{M}_n(\theta) \pas \mathcal{M}_{\cal E}(\theta), \;\;\;\;n
\rightarrow \infty, \eeq where $\mathcal{M}_n$ is the support
function of $Z_n = (b(\nu_n))^{-1}W_n.$ As
$$
\mathcal{M}_n(\theta) = \frac{1}{b(\nu_n)}\max_{\bk \in
T_n}\left\langle X_\bk,\theta\right\rangle,
$$
and
$$
E\left\langle X_\bk,\theta\right\rangle\left\langle
X_{\bll},\theta\right\rangle \rightarrow 0\;\;\; \mbox{when}\;\;
|\bk-\bll| \rightarrow \infty,
$$
we get (\ref{proofTh1})  by Lemma 12.
\halmos
\vspace{10pt}

{\bf Proof of Theorem \ref{thm2}.} For $h >0$, let us denote by
$C_{\bk,h}$ the cube\\ $[\bk h, (\bk+{\bf 1})h]$ and
$$G_n = \left\{\bk\,:\,C_{\bk,h}\cap T_n \neq \emptyset\right\},$$
$$T_{n,h} = \bigcup_{\bk \in G_n}C_{\bk,h}.$$
It is clear that $T_n \subset T_{n,h}$ and
$$
\l^m(T_{n,h}\setminus T_n) \leq \l^m((\partial T_n)^{2\sqrt{m}h}),
$$
and also
$$
\card \left\{G_n\right\}h^m = \l^m(T_{n,h}).
$$
It follows from (\ref{condVH}) that $\tilde{\nu}_n = \card
\left\{G_n\right\} \sim \nu_nh^{-m},$ therefore, $b(\tilde{\nu}_n)
\sim b(\nu_n).$

Let $Z_n = \left\{X_\bt, \bt\in T_n\right\},\;\;Z_{n,h} =
\left\{X_{\bk h}, \bk \in G_n \right\}.$ Then \beq\label{dist}
Z_{n,h} \subset (Z_n)^{d_n},\;\;\; Z_n \subset (Z_{n,h})^{d_n},
\eeq where $d_n = \max_{\bk \in G_n} \zeta_\bk,$ and $ \zeta_\bk
= \sup \left\{|X_\bt - X_{\bk h}|, \, \bt \in C_{\bk,h}\right\}. $

Relations (\ref{dist}) mean that
 \beq\label{dist1} \rho_{\bbb} (Z_n,\;Z_{n,h} ) \leq d_n.
\eeq By stationarity the random variables $\zeta_\bk $ are
identically distributed. From the continuity of $X$ it follows that
$\zeta_\bk < \infty$ a.s. As $\zeta_{\bk} $  is the supremum of
Gaussian random variables with variances less than $2\sigma^2(h),$
where
$$
\sigma^2(h)= \sup_{|\bt-\bs|\leq h}E|X_\bt - X_\bs|^2,
$$
then,
 according to the Fernique theorem from (\cite{Fernique}), for all $a <
\frac{1}{4\sigma^2(h)}$,
$$
M(a)= E\exp\{a\zeta_\bk^2\} <\infty.
$$
Now, due to Lemma \ref{propD}, we have a.s. \beq\label{dist2}
\limsup_n\{d_n\} \leq \sigma(h). \eeq By Theorem 1 for any $h>0$
\beq\label{dist3} \frac{1}{b(\nu_n)}Z_{n,h}\;\;\pas \cal E. \eeq
From (\ref{dist1}) and (\ref{dist2}) we have a.s.
$$
\limsup_n\rho_{\bbb}\left(\frac{1}{b(\nu_n)}Z_n,\;\;\frac{1}{b(\nu_n)}Z_{n,h}\right)\leq
\sigma(h).
$$
Hence for each $h$ a.s.
$$
\limsup_n\rho_{\bbb}\left(\frac{1}{b(\nu_n)}Z_n,\;\;\cal E\right)
\leq \sigma(h).
$$
Due to the  continuity we have that $ \sigma(h) \rightarrow 0,$ if
$h\rightarrow 0,$ therefore, finally we get
$$
\limsup_n\rho_{\bbb}\left(\frac{1}{b(\nu_n)}Z_n,\;\;\cal
E\right)\pas 0.
$$
The theorem is proved.
  \halmos
\vspace{10pt}

\vspace{10pt}

{\bf Proof of Theorem \ref{esper}.}
Due to the continuity of $f$ and the convergence (\ref{convW}) the
result of Theorem \ref{esper} will follow from the uniform
integrability of the family $\left\{f\left(\frac{W_n}
{b(\nu_n)}\right)\right\}.$

Due to the condition (\ref{moment}) we have
$$
f\left(\frac{W_n}{b(\nu_n)}\right) \;\;\leq \;\; C\left(\frac{D_n}{b(\nu_n)}\right)^p,
$$
where $D_n = \max_{\bk,\bll\in T_n}\left\|X_\bk - X_\bll\right\| \leq 2 \max_{\bk\in T_n}\left\|X_\bk\right\|.
$
Hence it is sufficient to state that for all $a>0$
\begin{equation}
\label{supE}
\sup_n E\exp{\left\{a\left(\frac{D_n}{b(\nu_n)}\right)^2\right\}} \;\;< \infty.
\end{equation}
 The latter relation  follows directly from Lemma \ref{propD}, and the
theorem is proved.\halmos


\section*{References}

\bibliographystyle{plain}
%


\end{document}